\documentclass[leqno,a4paper]{article}

\usepackage{amsmath,amssymb,amsthm,xypic,here}
\usepackage[mathscr]{eucal}

\numberwithin{equation}{section}

\theoremstyle{plain}
\newtheorem{theorem}{Theorem}[section]
\newtheorem{lemma}{Lemma}[section]
\newtheorem{proposition}{Proposition}[section]
\newtheorem{corollary}{Corollary}[section]

\newtheorem{question}{Question}[section]

\theoremstyle{remark}
\newtheorem{remark}{Remark}[section]

\newcommand{\Z}{\mathbf{Z}}

\newcommand{\F}{\mathbf{F}}
\newcommand{\nil}{\mathbf{nil}}
\newcommand{\e}{\epsilon}
\newcommand{\p}{_p^\wedge}

\newcommand{\bracket}[1]{\langle #1\rangle}

\newcommand{\SU}{\mathrm{SU}}
\newcommand{\Sp}{\mathrm{Sp}}

\newcommand{\Spin}{\mathrm{Spin}}
\renewcommand{\P}{\mathscr{P}}

\title{Homotopy nilpotency in $p$-compact groups}

\author{Daisuke Kishimoto\\Dept. of Math., Kyoto Univ.\\606-8502, Kyoto, Japan\\
         \texttt{kishi@math.kyoto-u.ac.jp}
 \and    
Shizuo Kaji\\Dept. of Math., Kyoto Univ.\\606-8502, Kyoto, Japan\\
         \texttt{kaji@math.kyoto-u.ac.jp}}

\begin{document}

\maketitle

\begin{abstract}
A $p$-compact group is a mod $p$ homotopy theoretical analogue of a
 compact Lie group. It is determined the homotopy nilpotency class of a $p$-compact group having the homotopy
 type of the $p$-completion of the direct product of spheres.
\end{abstract}

\section{Introduction}

H-spaces have been of great interest in algebraic topology and their multiplicative structures have been extensively
studied. In particular, the study of the homotopy commutativity is very
successful in the line of Hubbuck's torus theorem \cite{Hub}: A homotopy commutative, connected finite H-space is
homotopy equivalent to a torus. Its mod $p$ analogue is obtained firstly
by Aguad\'e and Smith \cite{AS}, McGibbon \cite{M2}, Hemmi \cite{He},
Lin \cite{L} and others. On the one hand, it is important in studying
the homotopy commutativity to determine whether a given H-space, not
necessarily finite, is homotopy commutative or not. McGibbon \cite{M}
specializes to localized Lie groups and raises the following question in which the
$p$-localization is in the sense of Bousfield and Kan \cite{BK} and we
write it by $-_{(p)}$.

\begin{question}
For which primes $p$ is a $p$-localized Lie group
homotopy commutative with the canonical multiplication? 
\end{question}

McGibbon \cite{M} shows the following which completely answers the above
question when the Lie group is simply connected.

\begin{theorem}
[McGibbon \cite{M}]
\label{McGibbon}
Let $G$ be a compact, connected, simple Lie group of type
 $(n_1,\ldots,n_l)$ with $n_1\le\cdots\le n_l$.
\begin{enumerate}
\item $G_{(p)}$ is homotopy commutative if $p>2n_l$.
\item $G_{(p)}$ is not homotopy commutative if $p<2n_l$ except for $(G,p)=(\Sp(2),3)$, $(G_2,5)$.
\end{enumerate}
\end{theorem}

Note that we can ask the above question by replacing the
$p$-localization with the $p$-completion and  the same result as Theorem
\ref{McGibbon} holds since it only deals with $p$-torsion in the group
$[G\times G,G]$, where the $p$-completion is in the
sense of Bousfield and Kan \cite{BK} and we denote it by
$-\p$. A $p$-compact group is defined by Dwyer and Wilkerson \cite{DW} as a mod
$p$ homotopy theoretic analogue of a compact Lie group. Precisely, a
$p$-compact group is a $p$-complete loop space having the finite mod $p$
cohomology (See \cite{DW} for basic results on $p$-compact groups.). Then we can slightly generalize the above question also by
replacing a $p$-localized Lie group with a $p$-compact group (See \cite{Sa}). Now we ask a further question:

\begin{question}
How homotopy noncommutative is the $p$-compact group for a
 given prime $p$?
\end{question}

Once we can measure the homotopy noncommutativity, this question makes
sense and we do this by the homotopy nilpotency classes defined as
follows. Let $X$ be a grouplike space, that is, $X$ satisfies all the
axioms of groups up to homotopy, and let $\gamma:X\times X\to X$ denote
the commutator map of $X$. Define the map $\gamma_k:X^{k+1}\to X$ by
$$\gamma_k=\gamma\circ(1\times\gamma)\circ\cdots\circ(\underbrace{1\times\cdots\times
1}_{k-1}\times\gamma):X^{k+1}\to X,$$
where $X^n$ stands for the direct product of $n$ copies of $X$. The
homotopy nilpotency class of $X$, denoted $\nil X$, is the least integer
$n$ such that $\gamma_n$ is null homotopic. Namely, if $\nil X=n$, then
$X$ is a nilpotent group of class $n$. Note that we normalize $\nil X$
such that $\nil X=1$ if and only if $X$ is homotopy commutative. For basic
facts on homotopy nilpotent H-spaces, readers refer to \cite{Z}.

We call a $p$-compact group exotic if its Weyl group corresponds to the reflection
group in Clark-Ewing's list \cite{CE} for $p$ odd or $\mathrm{GL}_4(\F_2)$ for $p=2$ but not the
$p$-completion of a Lie group, where the $2$-compact group having
$\mathrm{GL}_4(\F_2)$ as its Weyl group is Dwyer and Wilkerson's
$2$-complete finite loop space $\mathrm{DI(4)}$. For $p$ odd, we write by $N(X)=n$ if the exotic
$p$-compact group $X$ corresponds to the reflection group of number $n$ in Clark-Ewing's list \cite{CE}.
Recently, Andersen, Grodal, M\o ller and Viruel \cite{AGMV} and Andersen
and Grodal \cite{AG} give a classification of $p$-compact groups as
follows. For a $p$-compact group $X$, there exist a Lie group $G$ and
the direct product of exotic $p$-compact groups $X'$ such that 
$$X\cong G\p\times X'.$$ 

The aim of this article is to determine the homotopy nilpotency class
of a connected $p$-compact group when $p$
is a regular prime, that is, when the $p$-compact group has the homotopy type of
the $p$-completion of the direct product of spheres. We say that a
grouplike space $X$ is of type $(n_1,\ldots,n_l)$ if 
$$X_{(0)}\simeq S^{2n_1-1}_{(0)}\times\cdots\times S^{2n_l-1}_{(0)},$$
where, by the Hopf theorem, this holds if $X$ has a finite rational homology.

Let $G$ be a compact, simply connected, simple Lie group of type $(n_1,\ldots,n_l)$
with $n_1\le\cdots\le n_l$. Then, for the classical result of Serre
\cite{S}, 
$$G\p\simeq(S^{2n_1-1})\p\times\cdots\times(S^{2n_l-1})\p$$ 
if and only if $p\ge n_l$. Moreover, by combining the result of Kumpel \cite{K}
and Wilkerson \cite{W} together with checking Clark-Ewing's list \cite{CE}, we can see that
$$X\simeq(S^{2n_1-1})\p\times\cdots\times(S^{2n_l-1})\p$$ 
if and only if $p>n_l$ for an exotic $p$-compact group $X$ of type $(n_1,\ldots,n_l)$ with
$n_1\le\cdots\le n_l$. 

If the $p$-completion of a compact, connected Lie group has the homotopy
type of the $p$-completion of the direct product of spheres, it is
equivalent to the $p$-completion of the direct product of a torus and
compact, simply connected, simple Lie groups as an H-space. Note that
$\nil(X_1\times\cdots\times X_n)=\max\{\nil X_1,\ldots,\nil X_n\}$ for
grouplike spaces $X_1,\ldots,X_n$. Then the homotopy nilpotency class
of a $p$-compact group having the homotopy type of the $p$-completion of
the direct product of spheres is completely determined by the following,
where Theorem \ref{Lie-group} holds if we replace the $p$-completion
with the $p$-localization as well as Theorem \ref{McGibbon}.

\begin{theorem}
\label{Lie-group}
Let $G$ be a compact, simply connected, simple Lie group of type
 $(n_1,\ldots,n_l)$ with $n_1\le\cdots\le n_l$. 
\begin{enumerate}
\item $\nil G\p=3$ if $n_l\le p\le\frac{3}{2}n_l$ for $(G,p)=(F_4,17)$, $(E_6,17)$, $(E_8,41)$, $(E_8,43)$ or
      $(\mathrm{rank}\,G,p)=(1,2)$.
\item $\nil G\p=2$ if $\frac{3}{2}n_l<p<2n_l$ or $(G,p)$ is in the above exceptional case.
\end{enumerate}
\end{theorem}

\begin{theorem}
\label{exotic}
Let $X$ be an exotic $p$-compact group of type $(n_1,\ldots,n_l)$ with
 $n_1\le\cdots\le n_l$ for $p$ odd.
\begin{enumerate}
\item $\nil X=1$ if $p>2n_l$.
\item $\nil X=2$ if $n_l<p<2n_l$ except for $N(X)={\rm 2b}$ and $(N(X),p)=(23,11),(30,31)$.
\item $\nil X=3$ if $(N(X),p)$ is in the above exceptional case.
\end{enumerate}
\end{theorem}

The organization of this article is as follows. In section 2, we will
decompose the iterated commutator $\gamma_k$ above and, by using the
classical result on the odd primary component of the homotopy groups of
spheres, we will  prove a result
analogous to Theorem \ref{Lie-group} and Theorem \ref{exotic} in a more general
setting but including some indeterminacy. In section 3, the above
indeterminacy for Lie groups will be fixed in case by case analysis. The
case of the classical groups will be an easy consequence of the result of
Bott \cite{B}. We will show a cohomological criterion for a Samelson
product being nontrivial. In the case of the exceptional Lie groups, by
calculating the action of $\P^1$, we will use this criterion to fix the
indeterminacy. 
In section 4, we also fix the indeterminacy of exotic $p$-compact groups
analogously to the exceptional
Lie groups using the above cohomological criterion.

\section{Commutator map and Samelson products}

We begin with the following easy lemma in which we normalize classes of
nilpotency groups such that a group is nilpotent of class one if and
only if it is abelian.

\begin{lemma}
\label{commutator-group}
Let $K$ be a group generated by $x_1,\ldots,x_l$. For each
 $y_1,\ldots,y_{k+1}\in K$, $[y_1,[\cdots[y_k,y_{k+1}]\cdots]]$ is presented
 by a product of conjugation of $[x_{i_1}^{\pm
 1},[\cdots$
$[x_{i_k}^{\pm 1},x_{i_{k+1}}^{\pm 1}]\cdots]]$ for $1\le
 i_1,\ldots,i_{k+1}\le l$, where $[-,-]$ stands for the commutator in $K$.
In particular, $K$ is nilpotent
 of class $<k$ if and only if
 $[x_{i_1}^{\pm 1},[\cdots[x_{i_k}^{\pm 1},x_{i_{k+1}}^{\pm 1}]$
$\cdots]]=e$ for each $1\le
 i_1,\ldots,i_{k+1}\le l$, where $e$ is unity of $K$.
\end{lemma}

\begin{proof}
Let $y_i=[x_{j_1^i}^{\nu_1^i},[\cdots[x_{j_{k-1}^i}^{\nu_{k-1}^i},x_{j_k^i}^{\nu_k^i}]\cdots]]$ and
 $z_i\in K$ for $\nu_1^i,\ldots,\nu_k^i=\pm 1$, $1\le j^i_1,\ldots,j_k^i\le l$ and $1\le i\le n$. Denote
 by $x^y$ the conjugation $yxy^{-1}$. By using the
 formula
$$[x,yz]=[x,y][x,z]^y,$$
one has
\begin{align*}
[w,\prod_{i=n}^ny_i^{z_i}]&=[w,y_1^{z_1}][w,\prod_{i=2}^ny_i^{z_i}]^{y_1^{z_1}}=[w,y_1^{z_1}][w,y_2^{z_2}]^{y_1^{z_1}}[w,\prod_{i=2}^ny_i^{z_i}]^{y_1^{z_1}y_2^{z_2}}=\cdots\\
&=\prod_{i=1}^n[w,y_i^{z_i}]^{\prod_{j=1}^{i-1}y_j^{z_j}}=\prod_{i=1}^n[w^{z_i^{-1}},y_i]^{z_i\prod_{j=1}^{i-1}y_j^{z_j}},
\end{align*}
where $z_0=y_0=e$.
Then Lemma \ref{commutator-group} follows from the formula
$$[uv,w]=[v,w]^u[u,w].$$
\end{proof}

Let $X$ be a group-like space. The (generalized) Samelson product of
$\alpha:A\to X$ and $\beta:B\to X$, denoted $\bracket{\alpha,\beta}$, is the composition $A\wedge
B\stackrel{\alpha\wedge\beta}{\longrightarrow}X\wedge X\stackrel{\bar{\gamma}}{\to}X$,
where $\bar{\gamma}$ is the reduced commutator map of $X$. The
commutator map and Samelson products are, of course, closely related. In
particular, we have:

\begin{proposition}
\label{nil-Samelson}
Let $X$ be a group-like CW-complex such that $X=X_1\times\cdots\times
X_l$ and let $\e_k:X_k\to X$ and $\pi_k:X\to X_k$ denote the
 inclusion and the projection respectively. Then $\nil X<k$ if
$$(\e_{i_1}\circ\pi_{i_1})^{\pm
 1}\circ\bracket{\e_{j_1}^{\pm 1},(\e_{i_2}\circ\pi_{i_2})^{\pm 1}\circ\bracket{\cdots(\e_{i_k}\circ\pi_{i_k})^{\pm
 1}\circ\bracket{\e_{j_k}^{\pm 1},\e_{j_{k+1}}^{\pm 1}}\cdots}}=0$$
for each $1\le i_1,\ldots,i_k,j_1,\ldots,j_{k+1}\le l$, where we write
the map $\alpha:A\to X$ followed by the inverse $X\to X$ by
 $\alpha^{-1}$.
\end{proposition}

\begin{proof}
Consider the group $[X^{k+1},X]$ on which the group
 structure is given by the pointwise multiplication. Denote the
 commutator in the group $[X^{k+1},X]$ by $[-,-]$. 
Then, by definition, $[\lambda_1,[\cdots[\lambda_k,\lambda_{k+1}]\cdots]]$ for
 $\lambda_1,\ldots,\lambda_{k+1}\in[X^{k+1},X]$ is the composition:
$$X^{k+1}\stackrel{\Delta}{\longrightarrow}X^{(k+1)^2}\xrightarrow{\lambda_1\times\cdots\times\lambda_{k+1}}X^{k+1}\stackrel{\gamma_k}{\longrightarrow}X,$$
where $\Delta:X^{k+1}\to X^{(k+1)^2}$ is the diagonal map and $\gamma_k=\gamma\circ(1\times\gamma)\circ\cdots\circ(\underbrace{1\times\cdots\times
 1}_{k-1}\times\gamma):X^{k+1}\to X$ is the iterated commutator map.
Let $\rho_i:X^{k+1}\to X$ denote the $i$-th projection.
Then it follows from
\begin{equation}
\label{diagonal}
(\rho_1\times\cdots\times\rho_{k+1})\circ\Delta=1_{X^{k+1}}
\end{equation}
that
\begin{equation}
\label{gamma_k}
\gamma_k=[\rho_1,[\cdots[\rho_k,\rho_{k+1}]\cdots]].
\end{equation}

Let us concentrate in the subgroup $\mathscr{K}$ of $[X^{k+1},X]$ generated by
 $\iota_i\circ\pi_i\circ\rho_j$ for $1\le i\le l$ and $1\le
 j\le k+1$. 
Note that we can assume 
$$1_X=(\e_1\circ\pi_1)\cdots(\e_n\circ\pi_n),$$
where the right hand side is given by the pointwise
 multiplication. Then, for \eqref{gamma_k}, we have $\gamma_k\in\mathscr{K}$ and hence, by applying Lemma \ref{commutator-group} to the
 group $\mathscr{K}$, we obtain that $\gamma_k=0$ if and only if
$$[(\e_{i_1}\circ\pi_{i_1}\circ\rho_{j_1})^{\pm
 1},[\cdots[(\e_{i_k}\circ\pi_{i_k}\circ\rho_{j_k})^{\pm
 1},(\e_{i_{k+1}}\circ\pi_{i_{k+1}}\circ\rho_{j_{k+1}})^{\pm
 1}]\cdots]]=0$$
for each $1\le i_1,\ldots,i_{k+1}\le l$ and $1\le j_1,\ldots,j_{k+1}\le
 k+1$. Thus, for \eqref{diagonal}, $\gamma_k=0$ if and only if
$$\gamma_k\circ(\e_{i_1}^{\pm 1}\times\cdots\times\e_{i_{k+1}}^{\pm 1})=0$$
for each $1\le i_1,\ldots,i_{k+1}\le l$. 
Since $X$ is a group-like CW-complex, the induced map
 $[\wedge^{k+1}X,X]\to[X^{k+1},X]$ from the pinching map $X^{k+1}\to\wedge^{k+1}X$ is monic (See Lemma 1.3.5 of
 \cite{Z}). Then it follows that $\gamma_k=0$ if and only if
$$\bracket{\e_{i_1}^{\pm 1},\bracket{\cdots\bracket{\e_{i_k}^{\pm
 1},\e_{i_{k+1}}^{\pm 1}}\cdots}}=0$$
for each $1\le i_1,\ldots,i_{k+1}\le l$. 

Analogously, one can deduce that $\bracket{\e_{i_1}^{\pm
 1},\bracket{\cdots\bracket{\e_{i_k}^{\pm 1},\e_{i_{k+1}}^{\pm
 1}}\cdots}}=0$ for each $1\le i_1,\ldots,i_{k+1}\le l$ if
$$(\e_{i_1}\circ\pi_i)^{\pm
 1}\circ\bracket{\e_{j_1}^{\pm 1},(\e_{i_2}\circ\pi_{i_2})^{\pm 1}\circ\bracket{\cdots(\e_{i_k}\circ\pi_{i_k})^{\pm
 1}\circ\bracket{\e_{j_k}^{\pm 1},\e_{j_{k+1}}^{\pm 1}}\cdots}}=0$$
for each $1\le i_1,\ldots,i_k,j_1,\ldots,j_{k+1}\le l$ and the proof is completed.
\end{proof}

\begin{remark}
Suppose that $X$ is a loop space and each $X_i$ is the $p$-localized or
 $p$-completed odd sphere. Then, by using the retraction of McGibbon
 \cite{M}, it is immediate only to see that $\gamma_k=0$ if and only if $\bracket{\e_{i_1},\bracket{\cdots\bracket{\e_{i_k},\e_{i_{k+1}}}\cdots}}=0$ for each $1\le i_1,\ldots,i_{k+1}\le l$.
\end{remark}

Let $p$ be an odd prime. Before specifying Proposition \ref{nil-Samelson} to our case, let us
recall from \cite{T} basic facts on the $p$-primary component of the
homotopy groups of spheres. Denote the $p$-primary component of a finite
group $G$ by ${}^pG$. Then we have
\begin{equation}
\label{homotopy-sphere}
{}^p\pi_{2n-1+i}(S^{2n-1})\cong\begin{cases}
\Z/p&i=2p-3\\
0&i\le 4p-7,i\ne 2p-3.
\end{cases}
\end{equation}
Let $\alpha_1(3)$ denote a generator of ${}^p\pi_{2p}(S^3)\cong\Z/p$ and
let $\alpha_1(n)=\Sigma^{n-3}\alpha_1(3)$. Then
\begin{equation}
\label{alpha_1}
\alpha_1(2n-1)\text{ is a generator of }{}^p\pi_{2n+2p-4}(S^{2n-1})\cong\Z/p.
\end{equation} 
Moreover, one has
\begin{equation}
\label{alpha_1^2}
\alpha_1(3)\circ\alpha_1(2p)\ne 0\text{ and }\Sigma^2(\alpha_1(3)\circ\alpha_1(2p))=0.
\end{equation}

Let $X$ be a $p$-complete, simply connected, group-like space of type
$(n_1,\ldots,n_l)$ with $n_1\le\cdots\le n_l$. Kumpel \cite{K} shows that if
$p>n_l-n_1+2$, then
$$X\simeq(S^{2n_1-1})\p\times\cdots\times(S^{2n_l-1})\p.$$ 
Put $p>n_l-\frac{n_1}{2}+1$. Denote the inclusion $(S^{2n_i-1})\p\to X$ and the
projection $X\to(S^{2n_i-1})\p$ by $\e_i$ and $\pi_i$ respectively.
Then, by \eqref{homotopy-sphere} and \eqref{alpha_1}, we have
\begin{equation}
\label{<-,->}
\pi_s\circ\bracket{\e_i,\e_j}=\begin{cases}
a\alpha_1(2n_s-1)&\text{if }n_i+n_j=n_s+p-1\\
0&\text{otherwise,}
\end{cases}
\end{equation}
here $a\in\Z/p$ is possibly $0$. Hence, for \eqref{alpha_1^2}, 
$\pi_s\circ\bracket{\e_i,(\e_t\circ\pi_t)\circ\bracket{\e_j,\e_k}}\ne
0$ if and only if 
\begin{equation}
\label{<-,<-,->>}
\begin{cases}
n_s=2,\;n_i+n_t=p+1,\;n_j+n_k=n_t+p-1\\
\pi_2\circ\bracket{\e_i,\e_t}\ne 0,\;\pi_t\circ\bracket{\e_j,\e_k}\ne 0.
\end{cases}
\end{equation}
Therefore, by Lemma \ref{nil-Samelson}, an easy inspection shows:

\begin{theorem}
\label{nil-grouplike}
Let $X$ be a $p$-complete, simply connected, group-like space of type
 $(n_1,\ldots,n_l)$ with $n_1\le\cdots\le n_l$.
\begin{enumerate}
\item $\nil X=1$ if $p>2n_l$.
\item $\nil X=1$ or $2$ if:
\begin{enumerate}
\item $\frac{3}{2}n_l<p<2n_l$.
\item $n_l-\frac{n_1}{2}+1<p\le 2n_l$ and $X$ does not satisfy \eqref{<-,<-,->>}.
\end{enumerate}
\item $\nil X=3$ if $n_1=2,n_l<p\le\frac{3}{2}n_l$ and $X$ satisfy \eqref{<-,<-,->>}.
\end{enumerate}
\end{theorem}

The classical result of James and Thomas \cite{JT} shows that if the
above group-like space $X$ is a loop space and $n_l-n_1+2<p<2n_l$ then
$X$ is not homotopy commutative. Thus we have:

\begin{corollary}
\label{nil-loop}
Let $X$ be a simply connected $p$-compact group of type $(n_1,\ldots,n_l)$ with
 $n_1\le\cdots\le n_l$.
\begin{enumerate}
\item $\nil X=1$ if $p>2n_l$.
\item $\nil X=2$ if:
\begin{enumerate}
\item $\frac{3}{2}n_l<p<2n_l$.
\item $n_l-\frac{n_1}{2}+1<p\le 2n_l$ and $X$ does not satisfy \eqref{<-,<-,->>}.
\end{enumerate}
\item $\nil X=3$ if $n_1=2,n_l<p\le\frac{3}{2}n_l$ and $X$ satisfy \eqref{<-,<-,->>}.
\end{enumerate}
\end{corollary}

In most cases, the above corollary reduces the proof of Theorem \ref{Lie-group} and Theorem
\ref{exotic} to examining \eqref{<-,<-,->>} case by case when
$n_l<p<\frac{3}{2}n_l$. 

Hereafter we will use the following notation. Let $X$ be
a compact group such that 
$$X\simeq(S^{2n_1-1})\p\times\cdots\times(S^{2n_l-1})\p.$$
Then we denote the inclusion $(S^{2n_i-1})\p\to X$ and the projection
$X\to(S^{2n_i-1})\p$ by $\e_i$ and $\pi_i$ respectively.

\section{Lie groups}

Theorem \ref{Lie-group} {\it 2} immediately follows from a dimensional
reason of Corollary
\ref{nil-loop} (See \cite{MT} for the types of simple Lie groups.). In the case
$(G,p)=(G_2,7),(F_4,13),$
$(E_6,13),(E_7,19),(E_8,31)$, Hamanaka and Kono
\cite{HK} show that $\pi_2\circ\bracket{\e_2,\e_{n}}\ne 0$ and then,
for Corollary \ref{nil-loop}, we have $\nil G\p=3$, where $n$ is the
largest entry in the type of $G$. The remaining cases are listed in the
following table and we will check \eqref{<-,<-,->>} in these cases,
where $[x]$ in the table stands for the largest integer less than or equal to $x$.
\renewcommand{\arraystretch}{1.2}
\begin{table}[H]
\centering
\begin{tabular}{cc}
group&prime\\\hline
$\SU(n)$&$n\le p\le\frac{3}{2}n$\\
$\Sp(n)$&$2n<p<3n$\\
$\Spin(n)$&$2[\frac{n}{2}]<p<3[\frac{n}{2}]$\\
$E_7$&$23$\\
$E_8$&$37$
\end{tabular}
\end{table}

\subsection{Classical groups}

\subsubsection{$\SU(n)$}

Let $p$ be a prime such that $p\ge n$. Then we have
$$\SU(n)\p\simeq(S^3)\p\times(S^5)\p\times\cdots\times(S^{2n-1})\p.$$
The classical result of Bott \cite{B} shows that if $i+j>n$, the order of
the Samelson product $\bracket{\e_i,\e_j}$ is a non-zero multiple of
\begin{equation}
\label{Bott}
\nu_p\left(\frac{(i+j-1)!}{(i-1)!(j-1)!}\right),
\end{equation}
where $\nu_p(p^kq)=p^k$ for $(p,q)=1$.

It follows from \eqref{homotopy-sphere} and \eqref{Bott} that
$\bracket{1_{\SU(2)^\wedge_2},1_{\SU(2)^\wedge_2}}\ne 0$ and then
$\nil\SU(2)^\wedge_2\ge 2$. Since $\SU(2)\cong S^3$,
$\pi_9(SU(2))\cong\Z/9$ (See, for example, \cite{T}.). Then
$\bracket{1_{\SU(2)^\wedge_2},\bracket{1_{\SU(2)^\wedge_2},1_{\SU(2)^\wedge_2}}}=0$
and therefore $\nil\SU(2)^\wedge_2=2$.

Put $2<p<n\le\frac{3}{2}p$. Then it also follows from
\eqref{homotopy-sphere} and \eqref{Bott} that
$\pi_2\circ\bracket{\e_n,\e_{p-n+1}}\ne
0$ and $\pi_{p-n+1}\circ\bracket{\e_n,\e_{2p-2n}}\ne 0$,
and hence, for \eqref{<-,<-,->>},
$\pi_2\circ\bracket{\e_n,(\e_{p-n+1}\circ\pi_{p-n+1})\circ\bracket{\e_n,\e_{2p-2n}}}\ne
0$. Therefore, by Lemma \ref{nil-Samelson} and Corollary
\ref{nil-loop}, we have obtained $\nil\SU(n)\p=3$.

Put $n=p$. This is the only one case which is not covered by Corollary
\ref{nil-loop}. By an analogous calculation to the above case, one has
$\pi_2\circ\bracket{\e_{p-1},(\e_2\circ\pi_2)\circ\bracket{\e_{p-1},\e_2}}\ne 0$ and then
$\nil\SU(p)\p\ge 3$. Recall from \cite{T} that
$\Sigma^2:{}^p\pi_{2n+2k}(S^{2n-1})\to{}^p\pi_{2n+2k+2}(S^{2n+1})$ is
the zero map. Then $\pi_s\circ\bracket{\e_i^{\pm},(\e_t\circ\pi_t)^{\pm
1}\circ\bracket{\e_j^{\pm 1},(\e_u\circ\pi_u)^{\pm
1}\circ\bracket{\e_k^{\pm 1},\e_l^{\pm 1}}}}=0$ for each $1\le
i,j,k,l,s,t,u\le p$ and hence it follows from Lemma \ref{nil-Samelson}
that $\nil\SU(p)\p=3$.

\subsubsection{$\Sp(n)$}

Let $p>2n$. Then we have
$$\Sp(n)\p\simeq(S^3)\p\times(S^7)\p\times\cdots\times(S^{4n-1})\p.$$
The result of Bott \cite{B} also shows that if $i+j>n$, the order of the
Samelson product $\bracket{\e_{2i},\e_{2j}}$ is a non-zero multiple of
$$\nu_p\left(\frac{(2i+2j-1)!}{(2i-1)!(2j-1)!}\right).$$
Then, by an analogous calculation to the case of $\SU(n)$, we can deduce
$\nil\Sp(n)\p=3$ if $2n<p<3n$.

\subsubsection{$\Spin(n)$}

Friedlander \cite{F} shows that if $p$ is an odd prime, there
exists an equivalence of H-spaces
$$\Spin(2n+1)\p\simeq\Sp(n)\p.$$
Then it follows from the above result on $\Sp(n)$ that
$\nil\Spin(2n+1)\p=3$ if $2n<p<3n$. Harris \cite{H} shows that the
fibration sequence
$$\Spin(2n)\p\to\Spin(2n+1)\p\to(S^{2n+1})\p$$
splits and thus, by Corollary \ref{nil-loop}, $\nil\Spin(2n)\p=3$ if
$2n<p<3n$.

\subsection{Exceptional groups}

We first show a way to find a non-trivial Samelson product by the
Steenrod operation which is used in \cite{KO} and also in \cite{HK}. Let $X$ be a $p$-compact group such that
$$X=(S^{2n_1-1})\p\times\cdots\times(S^{2n_l-1})\p.$$
Then it is obvious that
$$H^*(BX;\F_p)=\F_p[x_{n_1},\ldots,x_{n_l}],\;|x_i|=2i.$$

\begin{lemma}
\label{cohomology-Samelson}
Let $\theta$ be the Steenrod operation. If $\theta
 x_i=ax_jx_k+$other terms with $a\ne 0$, then
 $\bracket{\e_j,\e_k}\ne 0$.
\end{lemma}

\begin{proof}
Suppose that $\bracket{\e_j,\e_k}=0$. Let
 $\mathrm{ad}:[V,\Omega W]\to[\Sigma V,W]$ denote the adjoint
 isomorphism. Then the Whitehead product
 $[\mathrm{ad}\e_j,\mathrm{ad}\e_k]=\mathrm{ad}\bracket{\e_j,\e_k}=0$
 and hence there exists a map $f:(S^{2n_j})\p\times(S^{2n_k})\p\to BX$
 such that the following square diagram is homotopy commutative.
$$\xymatrix{(S^{2n_j})\p\vee(S^{2n_k})\p\ar[rr]^{\mathrm{ad}\e_j\vee\mathrm{ad}\e_k}\ar[d]_\iota&&BX\vee
 BX\ar[d]^\nabla\\
(S^{2n_j})\p\times(S^{2n_k})\p\ar[rr]^f&&BX,}$$
where $\iota$ and $\nabla$ are the inclusion and the folding map. Note
 that we can assume
$$(\mathrm{ad}\e_m)^*(x_m)=u_{2n_m}$$
for $m=1,\ldots,l$, where $\sigma$ and $u_{n}$ denote the cohomology suspension and a
 generator of $H^n(S^n;\F_p)$ respectively. Then if $\theta
 x_i=ax_jx_k+$other terms with $a\ne 0$, 
$$f^*(\theta x_i)=au_{2n_j}\times u_{2n_k}\ne 0.$$
On the other hand, 
$$f^*(\theta x_i)=\theta f^*(x_i)=0$$
and this is a contradiction. Then Lemma \ref{P^1} is established.
\end{proof}

\begin{corollary}
\label{P^1}
Let $n_l<p\le\frac{3}{2}n_l$ and let $n_1<\cdots<n_l$. Suppose that
 $n_1=2$, $\P^1x_1=x_ix_l+$other terms and $\P^1x_l=x_jx_k+$other
 terms. Then $\nil X=3$.
\end{corollary}

\begin{proof}
For $n_1<\ldots<n_l$ and \eqref{homotopy-sphere}, one has
 $\pi_n\circ\bracket{\e_i,\e_l}=0$ unless $n=1$. On the other
 hand, it follows from Lemma \ref{cohomology-Samelson} that
 $\bracket{\e_i,\e_l}\ne 0$ and then
 $\pi_1\circ\bracket{\e_i,\e_l}\ne 0$. Analogously, one can see
 that $\pi_l\circ\bracket{\e_j,\e_k}\ne 0$ and then, for \eqref{<-,<-,->>} and Lemma
 \ref{nil-loop}, we have established Corollary \ref{P^1}.
\end{proof}

We next prepare notation for some symmetric polynomials. Denote $c_k$
the $k$-th elementary symmetric function in $t_1,\ldots,t_n$ for $k=1,\ldots,n$, that is,
$$1+c_1+\cdots+c_n=\prod_{i=1}^n(1+t_i).$$
Define symmetric polynomials $p_k$ and $s_k$ for $k=1,\ldots,n$ by
$$1-p_1+\cdots+(-1)^np_n=\prod_{i=1}^n(1-t_i^2)$$
and
$$s_k=t_1^{2k}+\cdots+t_n^{2k}$$
respectively.
Then one has Girard's formula
\begin{equation}
\label{Girard}
s_k=(-1)^kk\sum_{i_1+2i_2+\cdots+ni_n=k}(-1)^{i_1+\cdots+i_n}\frac{(i_1+\cdots+i_n-1)!}{i_1!\cdots
i_n!}p_1^{i_1}\cdots p_n^{i_n}
\end{equation}
(see \cite{MS}). In the canonical way, we will identify $c_k$ and $p_k$ with the universal
$k$-th Chern class and Pontrjagin class respectively.

\subsubsection{$E_7$ with $p=23$}

Recall that we have the commutative diagram:
\begin{equation}
\label{Spin(10)-E_6-E_7}
\xymatrix{
\Spin(10)\ar@{=}[r]\ar[d]_{i_1}&\Spin(10)\ar[d]^{i_2}\\
E_6\ar[r]^j&E_7
}
\end{equation}
Recall also that the mod $23$ cohomology of $B\Spin(10)$, $BE_6$ and $BE_7$ are:
\begin{align*}
H^*(B\Spin(10);\F_{23})&=\F_{23}[p_1,p_2,p_3,p_4,c_5],\\
H^*(BE_6;\F_{23})&=\F_{23}[x_2,x_5,x_6,x_8,x_9,x_{12}],\;|x_i|=2i,\\
H^*(BE_7;\F_{23})&=\F_{23}[y_2,y_6,y_8,y_{10},y_{12},y_{14},y_{18}],\;|y_i|=2i.
\end{align*}
In \cite{HK}, it is shown that $x_i$ and $y_i$ can be chosen such that:
\begin{align*}
j^*(y_i)&=x_i\;(i=2,6,8),&j^*(y_{10})&=x_5^2&j^*(y_{14})&=x_5x_9\\
i_1^*(x_2)&=p_1,&i_1^*(x_5)&=c_5&i_1^*(x_6)&=-6p_3+p_1p_2\\
i_1^*(x_8)&=12p_4+p_2^2-\tfrac{1}{2}p_1^2p_2
\end{align*}
Then, for \eqref{Spin(10)-E_6-E_7}, one has:
\begin{equation}
\label{Spin(10)-1}
\begin{aligned}
i_2^*(y_2)&=p_1&i_2^*(y_6)&=-6p_3+p_1p_2\\
i_2^*(y_8)&=12p_4+p_2^2-\tfrac{1}{2}p_1^2p_2&i_2^*(y_{10})&=c_5^4
\end{aligned}
\end{equation}
For a dimensional reason, one also has $i_1^*(x_9)=ap_1^2c_5-bp_2c_5$ for $a,b\in\F_{23}$
and then it follows from \eqref{Spin(10)-E_6-E_7} and \eqref{Spin(10)-1} that
\begin{equation}
\label{Spin(10)-2}
i_2^*(y_{14})=ap_1^2c_5^2-bp_2c_5^2.
\end{equation}

We can see from the following proposition that $(E_7)^\wedge_{23}$
satisfies the condition of Corollary \ref{P^1} and then $\nil(E_7)^\wedge_{23}=3$.

\begin{proposition}
$\P^1y_2=cy_6y_{18}+$other terms for $c\ne 0$ and
$\P^1y_6=dy_{10}y_{18}+ey_{14}y_{14}$+other terms for $d\ne 0$ or
 $e\ne 0$.
\end{proposition}

\begin{proof}
Define a ring homomorphism
$$\pi:\F_{23}[p_1,\ldots,p_4,c_5]\to\F_{23}[a_2,\ldots,a_4,b_5]/(a_2^3,a_3^2,a_4^2,b_5^3,12a_4+a_2^2)$$
 by
$$\pi(p_1)=0,\;\pi(p_i)=a_i\;(i=2,3,4),\;\pi(c_5)=b_5.$$
Then, for (\ref{Spin(10)-1}) and (\ref{Spin(10)-2}), one has
\begin{equation}
\label{pi-E_7}
\pi(i_2^*(y_2))=\pi(i_2^*(y_8))=\pi(i_2^*(y_6^2))=\pi(i_2^*(y_{14}y_{10}))=0.
\end{equation}
Put $\P^1y_4=cy_6y_{18}+$other terms for $c\in\F_{23}$. Then, for \eqref{pi-E_7} and a degree reason, one has
$$\pi(i_2^*(\P^1y_2))=c\pi(i_2^*(y_6y_{18})).$$
On the other hand, since $p_1=s_1$ and $\P^1s_1=2s_{12}$, Girard's formula (\ref{Girard}) yields that
$$\pi(i_2^*(\P^1y_2))=\pi(\P^1i_2^*(y_2))=\pi(\P^1p_1)=-15a_3a_4b_5^2\ne 0.$$
Thus we have $c\ne 0$.

We define a ring homomorphism
$$\pi':\F_{23}[p_1,\ldots,p_4,c_5]\to\F_{23}[a_2',a_4',b_5']/((a_2')^3,(a_4')^2,(b_5')^5,12a_4'+(a_2')^2)$$ 
by
$$\pi'(p_i)=0\;(i=1,3),\pi'(p_j)=a_j'\;(j=2,4),\;\pi'(c_5)=b_5'.$$
Then, for (\ref{Spin(10)-1}) and (\ref{Spin(10)-2}), we have
\begin{equation}
\label{pi'-E_7}
\pi'(i_2^*(y_2))=\pi'(i_2^*(y_6))=\pi'(i_2^*(y_8))=0.
\end{equation}
Put $\P^1y_6=dy_{10}y_{18}+ey_{14}y_{14}+$other terms for $d,e\in\F_{23}$. Then
it follows from \eqref{pi'-E_7} and a degree reason that
$$\pi'(i_2^*(\P^1y_6))=d\pi'(i_2^*(y_{10}y_{18}))+e\pi'(i_2^*(y_{14}y_{14})).$$
For Girard's formula (\ref{Girard}), we have:
\begin{align*}
\pi'(\P^1p_1)&=\pi'(\P^1s_1)=\pi'(2s_{12})=-a_2'(b_5')^4\\
\pi'(\P^1s_3)&=\pi'(6s_{14})=-9a_4'(b_5')^4
\end{align*}
Since $s_3=p_1^3-3p_1p_2+3p_3$, one has $\pi(\P^1p_3)=9a_4'(b_5')^4.$ and
 then, for (\ref{Spin(10)-1}), 
$$\pi'(\P^1(i_2^*(y_6)))=\pi'(\P^1(-6p_3+p_1p_2))=-19a_4'(b_5')^4\ne 0.$$
Therefore we have obtained $d\ne 0$ or $e\ne 0$.
\end{proof}

\subsubsection{$E_8$ with $p=37$}

Recall that the mod $37$ cohomology of $BE_8$ is given by
$$H^*(BE_8;\F_{37})=\F_{37}[z_2,z_8,z_{12},z_{14},z_{18},z_{20},z_{24},z_{30}],\;|z_i|=2i.$$
In order to see the action of $\P^1$ on $H^*(BE_8;\F_{37})$, we shall
choose suitable $y_i$.
 Let $\alpha_i$ ($i=1,\ldots,8$) and
$\tilde{\alpha}$ be respectively the simple roots and the dominant root
of $E_8$ as indicated in the following extended Dynkin diagram of $E_8$
(See \cite{MT} for details.).
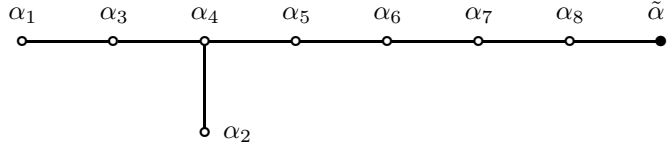
\begin{figure}[H]
\begin{center}
\setlength{\unitlength}{1.2cm}
\begin{picture}(7.5,1.7)
\thicklines
\multiput(0.5,1.1)(1,0){7}{\circle{0.1}}
\put(7.5,1.1){\circle*{0.124}}
\put(2.5,0.1){\circle{0.1}}
\multiput(0.55,1.1)(1,0){7}{\line(1,0){0.9}}
\put(2.5,1.05){\line(0,-1){0.9}}
\put(0.35,1.35){$\alpha_1$}
\put(1.35,1.35){$\alpha_3$}
\put(2.35,1.35){$\alpha_4$}
\put(3.35,1.35){$\alpha_5$}
\put(4.35,1.35){$\alpha_6$}
\put(5.35,1.35){$\alpha_7$}
\put(6.35,1.35){$\alpha_8$}
\put(7.35,1.35){$\tilde{\alpha}$}
\put(2.7,0.03){$\alpha_2$}
\end{picture}
\end{center}
\end{figure}
\noindent Denote the Weyl group of $E_8$ by $W$. We consider the
subgroup $K$ of $W$ which is generated by the reflections corresponding to
$\alpha_i$ ($i=2,\ldots,8$) and $\tilde{\alpha}$. Then, by choosing
appropriate generators $t_1,\ldots,t_8\in H^2(BT;\F_{23})$, one has
$$H^*(BT;\F_{37})^K=\F_{37}[p_1,\ldots,p_7,c_8],$$
where $T$ is a maximal torus of $E_8$. This is nothing but the
cohomology of $B\Spin(16)$ in $BE_8$.
Let $\varphi$ be the
elements of $W$ corresponding to $\alpha_1$. Then, by definition, $W$ is
generated by $K$ and $\varphi$. Hence one has
$$H^*(BT;\F_{37})^W=H^*(BT;\F_{37})^\varphi\cap\F_{37}[p_1,\ldots,p_7,c_8].$$
On the other hand, the canonical map $i:BT\to BE_8$ induces an isomorphism
$$i^*:H^*(BE_8;\F_{37})\stackrel{\cong}{\to}H^*(BT;\F_{37})^W.$$

It is shown in \cite{HK} that 
$$\varphi(c_1)=-c_1,\;\varphi(c_2)=c_2,\;\varphi(c_8)=c_8-\tfrac{1}{4}c_1c_7,\;\varphi(p_1)=p_1$$
and
\begin{equation}
\label{varphi}
\varphi(p_i)\equiv p_i+c_1h_i\mod (c_1^2)
\end{equation}
for $i=2,\ldots,7$, where
\begin{align*}
h_2&=\tfrac{3}{2}c_3&h_3&=-\tfrac{1}{2}(5c_5+c_2c_3)\\
h_4&=\tfrac{1}{2}(7c_7+3c_2c_5-c_3c_4)&h_5&=-\tfrac{1}{2}(5c_2c_7-3c_3c_6+c_4c_5)\\
h_6&=-\tfrac{1}{2}(5c_3c_8-3c_4c_7+c_5c_6)&h_7&=\tfrac{1}{2}(3c_5c_8-c_6c_7).
\end{align*}
Then it is immediate that
\begin{equation}
\label{z_2}
i^*(z_2)=p_1.
\end{equation}
By a direct calculation, Hamanaka and Kono \cite{H} show:

\begin{proposition}
[Hamanaka and Kono \cite{HK}]
\label{d_8d_12}
If $d_8\in H^{16}(BT;\F_{37})$ and $d_{12}\in
 H^{24}(BT;\F_{37})$ satisfy 
 $\varphi(d_8)\equiv d_8\mod{(c_1^2)}$ and
 $\varphi(d_{12})\equiv d_{12}\mod{(c_1^2,c_2^2)}$, then
$$d_8\equiv a\tilde{z}_8\mod(p_1^4),\;d_{12}\equiv b\tilde{z}_{12}\mod(p_1^2)$$
for $a,b\in\F_{37}$ and 
\begin{align*}
\tilde{z}_8&=120p_4+1680c_8+p_1^2p_2-36p_1p_3+10p_2^2,\\
\tilde{z}_{12}&=60p_6-p_1p_2p_3-5p_1p_5+\tfrac{5}{36}p_2^3-5p_2p_4+110p_2c_8+3p_3^2
\end{align*}
\end{proposition}

Since $z_2,z_8,z_{12}$ are algebraically independent, we have:

\begin{corollary}
[Hamanaka and Kono \cite{HK}]
\label{z_8z_12}
Let $\tilde{z}_8$ and $\tilde{z}_{12}$ be as in Proposition \ref{d_8d_12}.
Then we can choose generators $z_8$ and $z_{12}$ of $H^*(BE_8;\F_{37})$ such that
$$i^*(z_8)\equiv\tilde{z}_8,\;i^*(z_{12})\equiv\tilde{z}_{12}\mod(p_1^2).$$
\end{corollary}

Let us consider the generator $z_{14}$. For a dimensional reason, an element of degree $28$ in
$\F_{37}[p_1,\ldots,p_7,c_8]$ is given by a linear combination
$$ap_7+bp_2^2p_3+cp_2p_5+dp_3c_8+ep_3p_4\mod(p_1)$$
for $a,b,c,d,e\in\F_{37}$.
It is straightforward to check that
\begin{align*}
\varphi(p_2^2p_3)&\equiv p_2^2p_3+6c_1c_3^3c_4-12c_1c_3c_4c_6-10c_1c_4^2c_5,\\
\varphi(p_2p_5)&\equiv
 p_2p_5+3c_1c_3^2c_7+\tfrac{3}{2}c_1c_3c_5^2-c_1c_4^2c_5,\\
\varphi(p_3c_8)&\equiv p_3c_8-\tfrac{1}{4}c_1c_3^2c_7-\tfrac{5}{2}c_1c_5c_8+\tfrac{1}{2}c_1c_6c_7\\
\varphi(p_3p_4)&\equiv
 p_3p_4-\tfrac{1}{2}c_1c_3^3c_4+\tfrac{7}{2}c_1c_3^2c_7+c_1c_3c_4c_6\\
&\quad+5c_1c_3c_5^2-\tfrac{5}{2}c_1c_4^2c_5-5c_1c_5c_8-7c_1c_6c_7\mod{(c_1^2,c_2)}.
\end{align*}
Then it follows from \eqref{varphi} that:

\begin{proposition}
\label{d_14}
If $d_{14}\in H^{28}(BT;\F_{37})$ satisfy
 $\varphi(d_{14})\equiv d_{14}\mod{(c_1^2,c_2)}$, then
$$d_{14}\equiv a\tilde{z}_{14}\mod(p_1)$$
for $a\in\F_{37}$ and
$$\tilde{z}_{14}=480p_7-p_2^2p_3+40p_2p_5-12p_3p_4+312p_3c_8.$$
\end{proposition}

Since $z_2,z_8,z_{12},z_{14}$ are algebraically independent, we obtain:

\begin{corollary}
\label{z_14}
Let $\tilde{z}_{14}$ be as in Proposition \ref{d_14}. We can choose a generator
 $z_{14}$ of $H^*(BE_8;\F_{37})$ such that
$$i^*(z_{14})\equiv\tilde{z}_{14}\mod{(p_1)}.$$
\end{corollary}

We choose generators $z_2,z_8,z_{12},z_{14}$ of $H^*(BE_8;\F_{37})$ as
in \eqref{z_2}, Corollary \ref{z_8z_12} and Corollary \ref{z_14}. For
the following proposition, we can see that $(E_8)^\wedge_{37}$
satisfies the condition of Corollary \ref{P^1} and then $\nil(E_8)^\wedge_{37}=3$.

\begin{proposition}
$\P^1z_2=az_8z_{30}+$other terms for $a\ne 0$ and
$\P^1z_8=bz_{20}z_{24}+$other terms for $b\ne 0$.
\end{proposition}

\begin{proof}
Consider the ring homomorphism
$$\pi:\F_{37}[p_1,\ldots,p_7,c_8]\to\F_{37}[a_3,a_4,a_7,b_8]/(a_3^2,a_4^2,a_7^2,b_8^4,a_3a_4-26a_3b_8-40a_7)$$
defined by
$$\pi(p_i)=0\;(i=1,2,5,6),\;\pi(p_j)=a_j\;(j=3,4,7),\;\pi(c_8)=b_8.$$
Then, for \eqref{z_2}, Corollary \eqref{z_8z_12} and Corollary
 \ref{z_14}, we have
$$\pi(i_2^*(z_2))=\pi(i_2^*(z_{12}))=\pi(i_2^*(z_{14}))=0$$
and, for a degree reason, we also have
$$\pi(i_2^*(z_{18}))=0.$$
Put $\P^1z_2=az_8z_{30}+$other terms. Thus one has
$$\pi(i^*(\P^1z_2))=a\pi(i^*(z_8z_{30})).$$
On the other hand, it follows from Girard's formula (\ref{Girard}) and
 \eqref{z_2} that
$$\pi(i^*(\P^1z_2))=\pi(\P^1i^*(z_2))=\pi(\P^1p_1)=\pi(\P^1s_1)=\pi(2s_{19})=2a_4a_7b_8^2\ne 0$$
and then $a\ne 0$.

Define a ring homomorphism
$$\pi':\F_{37}[p_1,\ldots,p_7,c_8]\to\F_{37}[a_2',a_4',b_8']/((a_2')^2,(a_4')^6,(b_8')^6,a_4+14b_8)$$
by
$$\pi'(p_i)=0\;(i=1,3,5,6,7),\;\pi'(p_j)=a_j'\;(j=2,4),\;\pi'(c_8)=b_8'.$$
Then, for \eqref{z_2}, Corollary \ref{z_8z_12} and Corollary \eqref{z_14},  we have
$$\pi(i^*(z_4))=\pi(i^*(z_8))=\pi(i^*(z_{14}))=\pi(i^*(z_{12}^2))=0.$$
Put $\P^1z_8=by_{20}z_{24}+$other terms. Then we
 can see that
$$\pi(i_2^*(\P^1z_8))=b\pi(i^*(z_{20}z_{24})).$$
Let us make a direct calculation of $\pi(i_2^*(\P^1y_{16}))$. It follows
 from Girard's formula (\ref{Girard}) that:
\begin{align*}
\pi'(\P^1s_2)&=\pi'(4s_{20})=-6(b_8)^5&\pi'(\P^1s_4)=\pi'(8s_{22})=16a_2'(b_8')^5\\
\pi'(\P^1c_8)&=\pi'(s_{18}c_8)=26a_2'(b_8')^5
\end{align*}
Since $s_2=p_1^2-2p_2$ and $s_4=p_1^4-4p_1^2p_2+4p_1p_3+2p_2^2-4p_4$,
one has
$$\pi'(\P^1p_2)=3(b_8')^5,\;\pi'(\P^1p_4)=-a_2'(b_8')^5$$
and then, for Corollary \ref{z_8z_12}, 
\begin{align*}
\pi'(i^*(\P^1z_8))&=\pi'(\P^1i^*(z_8))\\
&=120\pi'(\P^1p_4)+1680\pi'(\P^1c_8)+20a_2'\pi'(\P^1p_2)=-3a_2'(b_8')^5\ne
 0.
\end{align*}
Hence we have obtained $b\ne 0$.
\end{proof}

\section{Exotic $p$-compact groups}

One can easily see that Theorem \ref{exotic} follows from
a dimensional reason of Corollary\ref{nil-loop} unless a $p$-compact group
corresponds to the reflection group in the following table.
\renewcommand{\arraystretch}{1.2}
\begin{table}[H]
\centering
\begin{tabular}{ccc}
number&degrees&prime\\\hline
2b&$2,n$&$n+1$\\
23&$2,6,10$&$11$\\
30&$2,12,20,30$&$31$
\end{tabular}
\end{table}
\noindent In the above table, the second column means half of the degrees of
generators of the invariant ring. Then it is the type of the
corresponding $p$-compact group.

Let $K$ be the reflection group in the above table and let
$\F_p[t_1,\ldots,t_l]^K=\F_p[x_1,\ldots,x_l]$, the invariant ring of
$K$, in which we put $|t_i|=2$. It follows from the above table that we
can set $|x_1|=4$. Note that we can
apply Corollary \ref{P^1} to the $p$-compact groups corresponding to
$K$ and the condition in Corollary \ref{P^1} is equivalent to:
\begin{equation}
\label{invariant-P^1}
\P^1x_1=x_ix_l+\text{other terms, }\P^1x_l=x_jx_k+\text{other terms}
\end{equation}
for some $i,j,k,l$, where the action of $\P^1$ comes from the relation $\P^1t_i=t_i^p$. Then we will show that $K$ satisfies this condition.

\subsection{$N(X)=$2b}

The reflection group of number 2b in Clark-Ewing's list \cite{CE} is
the Coxeter group $I_2(n)$ for $n\ge 3$ which corresponds to the
following Coxeter diagram and, in our case, $p=n+1$ (See \cite{Hu}.).
\begin{figure}[htbp]
\begin{center}
\setlength{\unitlength}{1.2cm}
\begin{picture}(1,0.3)
\thicklines
\multiput(0,0)(1,0){2}{\circle{0.1}}
\put(0.4,0.1){$n$}
\put(0.05,0){\line(1,0){0.9}}
\end{picture}
\end{center}
\end{figure}

\noindent Namely, it is the dihedral group of order $2n$, denoted $\mathscr{D}_{2n}$, acting on $\F_p\times\F_p$ via the matrices:
$$\begin{pmatrix}1&0\\0&-1\end{pmatrix},\;\begin{pmatrix}\cos\frac{2\pi}{n}&\sin\frac{2\pi}{n}\\\sin\frac{2\pi}{n}&-\cos\frac{2\pi}{n}\end{pmatrix},$$
here $\cos\frac{2\pi}{n}$ and $\sin\frac{2\pi}{n}$ exist in
$\F_p$ since $p=n+1$.
Consider the action of $\mathscr{D}_{2n}$ on $\F_p\times\F_p$ induced from the
above action via the linear transformation by
  $\begin{pmatrix}1&1\\-\sqrt{-1}&\sqrt{-1}\end{pmatrix}$. Then we obtain
$$\F_p[t_1,t_2]^{\mathscr{D}_{2n}}=\F_p[t_1t_2,t_1^n+t_2^n],$$
here $\sqrt{-1}$ also exists in $\F_p$.
Thus the invariant ring of $I_2(n)$ is
\begin{equation}
\label{I_2(n)-1}
\F_p[x_2,x_n],\;|x_i|=2i,\;\P^1x_2=x_2x_n
\end{equation}
and this satisfies \eqref{invariant-P^1}.

\subsection{$(N(X),p)=(23,11)$}

The reflection group of number $23$ in Clark-Ewing's list
\cite{CE} is the Coxeter group $H_3$ corresponding to the Coxeter
diagram: 
\begin{figure}[H]
\begin{center}
\setlength{\unitlength}{1.2cm}
\begin{picture}(2,0.3)
\thicklines
\multiput(0,0)(1,0){3}{\circle{0.1}}
\multiput(0.05,0)(1,0){2}{\line(1,0){0.9}}
\put(0.4,0.1){$5$}
\end{picture}
\end{center}
\end{figure}

\noindent By an analogous calculation of obtaining \eqref{I_2(n)-1}, we
have that the invariant ring of $I_2(5)$ over $\F_{11}$ is:
\begin{equation}
\label{I_2(n)-2}
P_1=\F_{11}[x_2,x_5],\;|x_i|=2i,\;\P^1x_2=x_2x_5^2-2x_2^6.
\end{equation}
Denote by
$P_2=\F_{11}[y_2,y_6,y_{10}],\;|y_i|=2i$ the invariant ring of $H_3$
over $\F_{11}$ in which we can set $i(y_2)=x_2$ by the map $i:P_2\to
P_1$ induced from the inclusion of the above Coxeter diagrams. 
Put $\P^1y_2=ay_2y_{10}$+other terms. Then it follows from
\eqref{I_2(n)-2} and a degree reason that
$$x_2x_5^2\equiv\P^1x_2\equiv\P^1i(y_2)\equiv i(\P^1y_2)\equiv ai(y_2)i(y_{10})\mod{(x_2^2)}$$
and hence $a\ne 0$. Therefore the invariant ring of $H_3$ over $\F_{11}$
satisfies \eqref{invariant-P^1}.

\subsection{$(N(X),p)=(30,31)$}

The reflection group of number 30 in
Clark-Ewing list \cite{CE} is the Coxeter group $H_4$ corresponding to
the Coxeter diagram:
\begin{figure}[H]
\begin{center}
\setlength{\unitlength}{1.2cm}
\begin{picture}(3,0.3)
\thicklines
\multiput(0,0)(1,0){4}{\circle{0.1}}
\multiput(0.05,0)(1,0){3}{\line(1,0){0.9}}
\put(0.4,0.1){$5$}
\end{picture}
\end{center}
\end{figure}
\noindent As above, we can see that the invariant ring of $I_2(5)$ over $\F_{31}$
is:
\begin{equation}
\label{I_2(n)-3}
Q_1=\F_{31}[x_2,x_5],\;|x_i|=2i,\;\P^1x_2\equiv x_2x_5^6,\;\P^1x_5\equiv
5x_5^7\mod{(x_2^2)}.
\end{equation}
Let $Q_2=\F_{31}[y_2,y_6,y_{10}],\;|y_i|=2i$ be the invariant ring of $H_3$ over
$\F_{31}$, where we put $i(y_2)=x_2$ by the canonical
by $i:Q_2\to Q_1$ as above. Set $\P^1y_2=ay_2y_{10}^3$+other
terms. Then, for \eqref{I_2(n)-3} and a degree reason, we have
$$x_2x_5^6\equiv\P^1x_2\equiv\P^1i(y_2)\equiv i(\P^1y_2)\equiv ai(y_2)i(y_{10}^3)
 \mod{(x_2^2)}$$
and hence $a\ne 0$ and we can put $i(y_{10})=x_5^2+$other terms. Thus we obtain
\begin{equation}
\label{y_10}
\P^1y_2=y_2y_{10}^3,\;\P^1y_{10}=10y_{10}^4\mod{(y_2^2,y_6)}.
\end{equation}
Denote by $Q_3=\F_{31}[z_2,z_{12},z_{20},z_{30}],\;|z_i|=2i$ the invariant
ring of $H_4$ over $\F_{31}$ in which we can set $j(z_2)=y_2$ by the canonical
map $j:Q_3\to Q_2$. Put
$\P^1z_2=bz_2z_{30}+cz_{12}z_{20}+$other terms. Then, for \eqref{y_10},
we have
\begin{equation}
\begin{aligned}
\label{y_2-y_10}
y_2y_{10}^3&\equiv\P^1y_2\equiv\P^1j(z_2)\equiv j(\P^1z_2)\\
&\equiv bj(z_2)j(z_{30})+cj(z_{12})j(z_{20})\mod{(y_2^2,y_6)}.
\end{aligned}
\end{equation}
Hence $b\ne 0$ or $c\ne 0$. Suppose that $b=0$. Then, for \eqref{y_2-y_10} we can set
$j(z_{20})=y_{10}^2$+other terms. Put $\P^1 z_{20}=dz_{20}z_{30}+$other
terms. Thus, for \eqref{I_2(n)-3} and a degree reason, we have
$$20y_{10}^5\equiv\P^1y_{10}^2\equiv\P^1j(z_{20})\equiv j(\P^1z_{20})\equiv
dj(z_{20})j(z_{30})\mod{(y_2^2,y_6)}$$
and hence $d\ne 0$.
Summarizing, we have established that
$\P^1z_2=bz_2z_{30}+cz_{12}z_{20}+$other terms with $b\ne 0$ or $c\ne 0$
such that if $b=0$, $\P^1z_{20}=dz_{20}z_{30}+$other terms with $d\ne
0$. This satisfies \eqref{invariant-P^1}.

\end{document}